\documentclass[11pt,reqno]{amsart}
\usepackage{amsthm}
\usepackage{amsmath}
\usepackage{amssymb}
\usepackage{amsbsy}
\usepackage{amsfonts}
\usepackage{latexsym}
\newtheorem{thm}{Theorem}[section]

\newtheorem{lem}[thm]{Lemma}

\theoremstyle{definition}
\newtheorem{defn}[thm]{Definition}
\theoremstyle{remark}
\newtheorem{rem}[thm]{Remark}
\numberwithin{equation}{section}
\makeatletter
\renewcommand{\@evenhead}{\thepage \hfil V.V.BORZOV \hfil}
\renewcommand{\@oddhead}{\hfil GENERALIZED HERMITE POLYNOMIALS  \hfil \thepage}
\makeatother
\begin{document}
\begin{center}

{\Large \bf GENERALIZED  HERMITE POLYNOMIALS
\footnote{This research was supported by RFFIgrant No 00-01-00500}}

\vspace{5mm}

{\large \bf V.V. Borzov}

\medskip

{\em Department of Mathematics, St.Petersburg University of
Telecommunications, \\ 191065, Moika  61,  St.Petersburg, Russia}

\medskip

%{\em (Received \qquad\qquad\qquad)}

\end{center}
\bigskip
The new method for obtaining a variety of extensions of Hermite
polynomials is given. As a first example  a family of orthogonal
polynomial systems which includes the generalized Hermite
polynomials is considered. Apparently,  either these polynomials
satisfy the differential equation of the second order obtained in
this work or there is no differential equation of a finite order
for these polynomials.

\bigskip


\medskip


\bigskip

\section{Introduction}

In our former work (\cite{bor}) we constructed an appropriate
oscillator algebra $A_\mu$ corresponding to the system of
polynomials which are orthonormal with respect to a measure $\mu$
in the space
\begin{math}{\tt H}_x=L^2(R^1;{\mu(dx))}\end{math}.
By a standard manner  the energy operator (hamiltonian)
$H_{\mu}={X_{\mu}^2}+{P_{\mu}^2}$ was defined. The position
operator $X_{\mu}$ was introduced by the recurrent relations of
the given polynomials system; the momentum operator $P_{\mu}$ was
determined as an unitary equivalent to the position operator
$Y_{\mu}$ in the dual space
\begin{math}{\tt H}_y=L^2(R^1;{\nu(dy))}\end{math}.
In (\cite{bor}) it was proved that the usual differential equations
for the classical polynomials are equivalent to the equations of
the form $H_{\mu}\psi_n=\lambda_n{\psi_n}$, where  the eigenvalues
of the corresponding hamiltonian $H_{\mu}$ denote by $\lambda_n$.
The central problem with a derivation of the  differential
 equations was finding a representation of the annihilation
operator ${a_{\mu}}^-$ (or another "reducing" operator)
of the algebra $A_\mu$ by a differential operator in the space
${\tt H}_x$. Unfortunately, these formulas (\cite{bord}) are
rather complicated in the general case. Therefore our interest is
in describing such orthogonal polynomials systems for which
appropriate representations are simple. On this basis one can
obtain some differential equations of a finite order (it is
desirable that we have to deal only with  differential equations
of the second order).

The results of this work may be thought of as a first step forward
in this direction. From our  point of view we consider a family of
orthogonal polynomial systems which includes the generalized
Hermite polynomials (\cite{sze}). These polynomials have been
studied extensively in the monograph (\cite{Chih}). Therefore the
polynomials of the considered family is called Hermite-Chihara
polynomials.

The paper is organized as follows. A generalized derivation
operator is introduced in Sec.2. By these operators a family of
the Hermite-Chihara polynomials is determined  in Sec.3. More
exactly the annihilation operator $a_{\mu}^{-}$ of the algebra
$A_\mu$ corresponding to the system of the polynomials may be
represented by a generalized derivation operator $D_{\vec{v}}$.
This operator is defined by a positive sequence  $\vec{v}$. In
what follows we shall call this  sequence  $\vec{v}$ a "governing
sequence". In Sec.4 we construct the generators $X_{\mu}$,
$P_{\mu}$ and $H_{\mu}$ of the algebra $A_\mu$ corresponding to a
system of the Hermite-Chihara polynomials. As an example of such
polynomials we consider the "classical"
Hermite-Chihara polynomials (\cite{Chih}) in Sec.5. Moreover, in
this section a new derivation of the well-known
(\cite{sze},\cite{Chih}) differential equation for these
polynomials is presented. Further, in Sec.6 we introduce a special
family of orthogonal systems of Hermite-Chihara polynomials which
includes the classical Hermite-Chihara polynomials. Furthermore,
in this section we construct a "governing sequence" $\vec{v}$ of
an appropriate generalized derivation operator $D_{\vec{v}}$. Then
in Sec.7 we obtain a differential equation of the second order for
above-mentioned polynomials by analogy with the derivation of the
differential equation given in Sec.5. Finally, in the conclusion
we consider the following conjecture. If the polynomials of a
system of orthogonal Hermite-Chihara polynomials satisfy a
differential equation of the second order, then these polynomials
belong to the special family of orthogonal systems of
Hermite-Chihara polynomials introduced in Sec.6. Moreover, the
other Hermite-Chihara polynomials do not satisfy any differential
equation of a finite order.

\section{Generalized derivation operator}

In this section we introduce a new class of differential operators
(they are the infinite order in general case) which play a large
role in the construction of the Hermite-Chihara polynomials. Let
$\vec{v}=\left\{ {v_{n}}\right\}_{n=0}^{\infty }$ be a monotone
nondecreasing sequence:
\begin{equation}
1=v_{0}\leq{v_{1}}\leq{v_{2}}\leq \dots\leq{v_{n}}\leq \dots.
\qquad \label{triv1}
\end{equation}
This sequence $\vec{v}$ define a linear operator $D_{\vec{v}}$ by
the relations:
\begin{equation}
{D_{\vec{v}}x^0}=0, \qquad {D_{\vec{v}}x^n}={v_{n-1}}x^{n-1},
\qquad n=1,2\dots ,
 \label{triv2}
\end{equation}
on the set of the formal power series of real argument $x$.

We will seek for the operator $D_{\vec{v}}$ of the type
\begin{equation}
D_{\vec{v}}=\sum_{n,m=0}^{\infty
}{a_{nm}{x^n}\frac{d^m}{dx^m}}\quad. \label{triv3}
\end{equation}
Substituting (\ref{triv3}) in (\ref{triv2}), we get the following
formula:
\begin{equation}
D_{\vec{v}}=\sum_{k=0}^{\infty
}{\varepsilon_{k}{x^{k-1}}\frac{d^k}{dx^k}}\quad, \label{triv4}
\end{equation}
The coefficients $\left\{ {\varepsilon_{k}}\right\}_{k=0}^{\infty
}$ are defined by the recurrent relations:
\begin{equation}
\varepsilon_1=v_0=1, \quad
\varepsilon_k=\frac{v_{k-1}}{k!}-\varepsilon_{k-1}-{\frac{\varepsilon_{k-2}}{2!}}-\cdots-
{\frac{\varepsilon_{1}}{(k-1)!}}, \quad k=1,2\dots . \label{triv5}
\end{equation}
\begin{defn}
A differential operator $D_{\vec{v}}$ determined by formulas
(\ref{triv4}),(\ref{triv5}) is called a generalized derivation
operator induced of the  sequence $\vec{v}$.
\end{defn}

\begin{lem}
For the order of a generalized derivation operator $D_{\vec{v}}$
defined by formulas (\ref{triv4}),(\ref{triv5}) to be finite it is
necessary and sufficient that the following equalities:
\begin{equation}
\varepsilon_{k+1}=\varepsilon_{k+2}=\cdots=0. \label{triv6}
\end{equation}
was valid.
\end{lem}

To take  three examples of  generalized derivation operators of a
finite order.

 1. Let $k=1$. There exist a unique solution of the
system (\ref{triv6}):
\begin{equation}
\vec{v}=\left\{n+1\right\}_{n=0}^{\infty }. \label{triv7}
\end{equation}
The generalized derivation operator $D_{\vec{v}}$  corresponding
to $\vec{v}$ take the following form:
\begin{equation}
D_{\vec{v}}=\frac{d}{dx}\quad. \label{triv8}
\end{equation}

2. Let $k=2$ and let $v_1$ to be a number such that  $v_1\geq1$.
There exist a one-parameter family of the solution
$\vec{v}=\left\{ {v_{n}}\right\}_{n=0}^{\infty }$ of the system
(\ref{triv6}):
\begin{equation}
v_0=1, \qquad v_n={C_{n+1}^2}{v_1}-n^2+1, \qquad n=1,2\dots .
\label{triv9}
\end{equation}
The  generalized derivation operator $D_{\vec{v}}$ corresponding
to $\vec{v}$ take the following form:
\begin{equation}
D_{\vec{v}}={\frac{d}{dx}}+x({\frac{v_1}{2}}-1){\frac{d^2}{dx^2}}\quad.
\label{triv10}
\end{equation}
If $v_1=4$, then from (\ref{triv9}),(\ref{triv10}) we have
\begin{equation}
\vec{v}=\left\{(n+1)^2\right\}_{n=0}^{\infty },\qquad
D_{\vec{v}}={\frac{d}{dx}}+x{\frac{d^2}{dx^2}}\quad.
\label{triv11}
\end{equation}

3. Let $k=3$ and let  ${v_1},{v_2}$ to be some number such that
$1\leq{v_1}\leq{v_2}$. There exist a two-parameter family of the
solution $\vec{v}=\left\{ {v_{n}}\right\}_{n=0}^{\infty }$ of the
system (\ref{triv6}):
\begin{align}
v_0=1\leq{v_1}, \qquad
v_n={C_{n+1}^3}{v_{2}}-{\frac{(n+1)n(n-2)}{2}}{v_1}+\nonumber\\
+{\frac{(n+1)(n-1)(n-2)}{2}},\qquad n=1,2\dots . \label{triv12}
\end{align}
The generalized derivation operator $D_{\vec{v}}$ corresponding to
$\vec{v}$ take the following form:
\begin{equation}
D_{\vec{v}}={\frac{d}{dx}}+x({\frac{v_1}{2}}-1){\frac{d^2}{dx^2}}+{x^2}
{\frac{v_2-3v_1+3}{3!}}{\frac{d^3}{dx^3}}\quad. \label{triv13}
\end{equation}
If  $v_1=8$,\quad$v_2=27$, then from (\ref{triv12}) and
(\ref{triv13}) we have
\begin{equation}
\vec{v}=\left\{(n+1)^3\right\}_{n=0}^{\infty },\qquad
D_{\vec{v}}={\frac{d}{dx}}+x{\frac{d^2}{dx^2}}+{x^2}
{\frac{d^3}{dx^3}}\quad. \label{triv14}
\end{equation}

\section{Hermite-Chihara polynomials}

Let $\mu$ be a symmetric probability measure, i.e.  the all odd
moments of the measure $\mu$ are vanish and
$\int_{-\infty}^{\infty}{\mu(dx)}=1$.
 In this section we consider a system of polynomials , which are
orthonormal  with  respect to the measure $\mu$ , such that there
is a representation of the annihilation operator of the
 oscillator algebra $A_\mu$ corresponding to this system by a generalized
derivation operator.

Recall (\cite{bor}) that the recurrent relations of a canonical
orthonormal polynomials system $\left\{
{\psi_{n}(x)}\right\}_{n=0}^{\infty }$ take the following form:
\begin{equation}
{x {\psi_{n}(x)}}=
{b_{n}{\psi_{n+1}(x)}}+{b_{n-1}{\psi_{n-1}(x)}},\qquad  n\geq 1,
\label{triv15}
\end{equation}
\begin{equation}
\psi_{0}(x)=1,\qquad \psi_{1}(x)=\frac{x}{b_0}. \label{triv16}
\end{equation}
In (\cite{bor})it was described how to get the positive sequence
$\left\{ {b_{n}}\right\}_{n=0}^{\infty }$ from the given sequence
$\left\{ {\mu_{2n}}\right\}_{n=0}^{\infty }$ of  even moments of a
symmetric positive measure $\mu$.

The question we are interested now is when for a canonical
orthonormal polynomials system $\left\{
{\psi_{n}(x)}\right\}_{n=0}^{\infty }$ there are two sequences
such that:

1. a positive sequence $\vec{v}=\left\{
{v_{n}}\right\}_{n=0}^{\infty }$ which satisfies (\ref{triv1});

2. a real sequence $\vec{\gamma}=\left\{
{\gamma_{n}}\right\}_{n=0}^{\infty }$ for which are hold the
following relations:
\begin{equation}
{D_{\vec{v}}\psi_0}=0, \qquad
{D_{\vec{v}}\psi_n}={\gamma_{n}}\psi_{n-1}, \qquad n=1,2\dots ,
\label{triv17}
\end{equation}
where the generalized derivation operator $D_{\vec{v}}$ is
determined by formulas (\ref{triv4}),(\ref{triv5}).

We denote by $[n]$ the following symbol:
\begin{equation}
[0]=0, \qquad [n]=\frac{b_{n-1}^2}{b_{0}^2}, \qquad n=1,2\dots .
\label{triv18}
\end{equation}
Let $J$ be a symmetric Jacobi matrix
 $$
J=\left\{ b_{ij} \right\} _{i,j=0}^{\infty }
 $$
which has the positive elements $b_{i,i+1}=b_{i+1,i},\quad
i=0,1,$\dots only distinct from zero. Then the polynomials of the
first kind can be represented in the form (\cite{bodak}):
\begin{equation}
\psi_{n}(x)={\sum_{m=0}^{\epsilon(\frac{n}{2})}}\frac{(-1)^{m}}{\sqrt{[n]!}}
b_{0}^{2m-n} {\alpha_{2m-1,n-1}}x^{n-2m}, \label{triv19}
\end{equation}
where the greatest integer function is denoted by
$\epsilon(\alpha)$  . The coefficients $\alpha_{2m-1,n-1}$ for any
$n\geq1,\quad \epsilon(\frac{n}{2})\geq{m}\geq1$ are defined by
the following equalities:
\begin{equation}
\alpha_{-1,n-1}=0,\qquad
\alpha_{2m-1,n-1}=\sum_{k_1=2m-1}^{n-1}[k_1]\sum_{k_2=2m-3}^{k_1-2}[k_2]\cdots
\sum_{k_m=1}^{k_{m-1}-2}[k_m]. \label{triv20}
\end{equation}
Substituting (\ref{triv20}),(\ref{triv19}) and (\ref{triv18}) into
(\ref{triv17}), it is easy to prove the following theorem
(\cite{bord}).
\begin{thm}
\label{a} Let the orthonormal polynomial system $\left\{
{\psi_{n}(x)}\right\}_{n=0}^{\infty }$ is defined by
(\ref{triv20}),(\ref{triv19}) and (\ref{triv18}). For existence
two sequences  $\vec{v}=\left\{ {v_{n}}\right\}_{n=0}^{\infty }$
and $\vec{\gamma}=\left\{ {\gamma_{n}}\right\}_{n=0}^{\infty }$
such that the conditions (\ref{triv17}) are hold it is necessary
and sufficient that

1. the sequence $\vec{v}=\left\{ {v_{n}}\right\}_{n=0}^{\infty }$
satisfies (\ref{triv1}) and the following conditions:
\begin{equation}
v_{n-2}v_{2p-1}+v_{2p-3}v_{n-2p}=v_{n}v_{2p-3}+v_{2p-1}v_{n-2p},
\label{triv21}
\end{equation}
for any $n\geq2,\quad 2p\leq{n};$

2. the coefficients  $\alpha_{2m-1,n-1}$ take the following form
\begin{equation}
\alpha_{2m-1,n-1}=\frac{[2m-1]!!(v_{n-1})!}{(v_{2m-1})!(v_{n-2m-1})!},\qquad
(v_k)!={v_0}{v_1}\cdots{v_k}, \label{triv22}
\end{equation}
as $n\geq1,\quad 2m\leq{n}$ and regarding $(v_{-1})!=(v_0)!=1.$

In this case the sequence  $\vec{\gamma}=\left\{
{\gamma_{n}}\right\}_{n=0}^{\infty }$ is defined by the following
formulas:
\begin{equation}
\gamma_n=\sqrt{\frac{{v_1}{v_{n-1}}}{b_0^2(v_n-v_{n-2})}},\qquad
n\geq1. \label{triv23}
\end{equation}
\end{thm}

Here we will not given the proof of this theorem (see \cite{bord})
to save room. However we present some formulas arising from the
proof. These expressions relate the sequence
 $\vec{\gamma}=\left\{
{\gamma_{n}}\right\}_{n=0}^{\infty }$ and the coefficients
$\alpha_{2m-1,n-1}$:
\begin{align}
\label{triv24}
\gamma_1=\frac{1}{b_0},\qquad{b_0\frac{\sqrt{[2]}}{2}\gamma_2}&=\varepsilon_1+\varepsilon_2,\\
\label{triv25}
 {\sqrt{[2p+1]}\gamma_{2p+1}}&={\frac{\alpha_{2p-1,2p}}{[2p-1]!!}}\gamma_1,\\
\label{triv26}
 {\sqrt{[2p+2]}\gamma_{2p+2}}&={\frac{\alpha_{2p-1,2p+1}}{\alpha_{2p-1,2p}}}\sqrt{[2]}\gamma_2,
\end{align}
where $\varepsilon_1$ and $\varepsilon_2$ are defined by
(\ref{triv5}).

Now we shall give the following definition.

\begin{defn}
The orthonormal polynomials system $\left\{
{\psi_{n}(x)}\right\}_{n=0}^{\infty }$ completed in ${\tt
H}_x={L^2}(R;\mu(dx))$ is called a system of Hermite-Chihara
polynomials if these polynomials are defined by
(\ref{triv20}),(\ref{triv19}) and (\ref{triv18}).
\end{defn}

\begin{rem}
It is clear that the Hermite polynomials fall in this class. Here
$$ \vec{v}=\left\{n+1\right\}_{n=0}^{\infty
},\quad{D_{\vec{v}}=\frac{d}{dx}},\quad{b_{n-1}^2}=\frac{n}{2},\quad[n]=n,\quad{n\geq1.}
 $$
According to (\ref{triv22}), we have
\begin{equation}
\alpha_{2m-1,n-1}=\frac{n!}{{2^m}m!(n-2m)!}.\qquad \label{triv27}
\end{equation}
Substituting (\ref{triv27}) into (\ref{triv19}), we obtain the
usual form of the Hermite polynomials (see \cite{bor}, \cite{sze}
). In addition, $\gamma_n=\sqrt{2n},\quad{n\geq1}$ and then
(\ref{triv17}) is reduced to the usual rule of derivation for the
 Hermite polynomials:
\begin{equation}
\frac{d}{dx}{H_n(x)}=2nH_{n-1}(x)\quad. \label{triv28}
\end{equation}
\end{rem}
\begin{rem}
According to the theorem {\ref{a}}, by any sequence
$\vec{v}=\left\{ {v_{n}}\right\}_{n=0}^{\infty }$ complying with
(\ref{triv1}) and (\ref{triv21}) we can write the the coefficients
$\alpha_{2m-1,n-1}$ which take the following form:
\begin{align}
\label{triv29}
[1]=1,\qquad[n]&=\frac{v_{n-1}(v_n-v_{n-2})}{v_1},\qquad{n\geq2},\\
\label{triv30}
 {b_{n-1}^2}&={b_0^2}\frac{v_{n-1}(v_n-v_{n-2})}{v_1},\qquad{n\geq2}.
\end{align}
The polynomials $\psi_{n}(x)$ satisfy the recurrent relations
(\ref{triv15}) and (\ref{triv16}). By solving the Hamburger moment
problem of the Jacobi matrix $J$, we obtain the symmetric
probability measure $\mu$ such that the polynomials of the system
$\left\{ {\psi_{n}(x)}\right\}_{n=0}^{\infty }$ are orthonormal
with respect to $\mu$. If the moment problem for the Jacobi matrix
$J$ is a determined one, then the measure $\mu$ is defined
uniquely. Otherwise (when the  moment problem for the Jacobi
matrix $J$ is a undetermined one) there is a infinite family of
such measures (see
 \cite{akh}).
\end{rem}

\section{Oscillator algebra for the
Hermite-Chihara polynomials}

In this section we construct the generalized Heisenberg algebra
$A_\mu$ corresponding to the system of the
Hermite-Chihara polynomials (see \cite{bor}).

Let $\vec{v}=\left\{ {v_{n}}\right\}_{n=0}^{\infty }$ be the
positive sequence such that the conditions (\ref{triv1}) and
(\ref{triv21}) are hold. Then the sequence $\left\{
{b_{n}}\right\}_{n=0}^{\infty }$ can be found by (\ref{triv30}).
Furthermore, we obtain the system of the Hermite-Chihara
polynomials $\left\{ {\psi_{n}(x)}\right\}_{n=0}^{\infty }$ by the
formulas (\ref{triv19}) and (\ref{triv22}). These polynomials
satisfy the recurrent relations (\ref{triv15}) and (\ref{triv16})
with above-mentioned coefficients $\left\{
{b_{n}}\right\}_{n=0}^{\infty }$. Under the condition
\begin{equation}
\sum_{n=0}^{\infty}{b_{n}^{-1}}=
\frac{v_1}{b_0}\sum_{n=1}^{\infty}{\frac{1}{\sqrt{v_{n}(v_{n+1}-v_{n-1})}}}
=\infty. \label{triv31}
\end{equation}
the moment problem for the corresponding Jacobi matrix is a
determined one (see \cite{akh}). There is the only symmetric
probability measure $\mu$ such that the polynomials  $\left\{
{\psi_{n}(x)}\right\}_{n=0}^{\infty }$ are orthonormal in the
space ${H}_x={L^2}(R;\mu(dx))$ . In addition, the even moments
$\mu_{2n}$ of the measure $\mu$ can be found from the following
algebraic equations system ($b_{-1}=0,\quad{n\geq0}$)
\begin{equation}
\sum_{m=0}^{\epsilon(\frac{n}{2})}{\sum_{m=0}^{\epsilon(\frac{n}{2})}{\frac{(-1)^{m+s}}{(b_{n-1})!}}}
{\alpha_{2m-1,n-1}}{\alpha_{2s-1,n-1}}{\mu_{2+2n-2m-2s}}={b_{n-1}^2}+{b_{n}^2}.
\label{triv32}
\end{equation}

It is easy to check that the condition (\ref{triv31}) is correct
for the "classical"  Hermite-Chihara polynomials to be considered
in the next section.

From theorem {\ref{a}} it follows that there are the generalized
derivation operator $D_{\vec{v}}$ determined  for given sequence
$\vec{v}$ by formulas (\ref{triv4}),(\ref{triv5}) and the sequence
$\vec{\gamma}=\left\{ {\gamma_{n}}\right\}_{n=0}^{\infty }$:
\begin{equation}
\gamma_n=\frac{v_{n-1}}{b_{n-1}},\qquad{n\geq1}, \label{triv33}
\end{equation}
such that the relations (\ref{triv17}) are valid .

Using the methods of (\cite{bor}), we construct the generalized
Heisenberg algebra $A_{\mu}$ corresponding to the  ortonormal
system
 $\left\{ {\psi_{n}(x)}\right\}_{n=0}^{\infty
}$. By the usual formulas we define ladder operators $a_{\mu}^{-}$
(the annihilation operator), $a_{\mu}^{+}$ (the creation operator)
 and the number operator $N$.
It is readily seen that
\begin{equation}
a_{\mu}^{-}=D_{\vec{v}}f(N), \label{triv34}
\end{equation}
 The operator-function $f(N)$ acts on basis vectors $\left\{
{\psi_{n}(x)}\right\}_{n=0}^{\infty }$ by formulas:
\begin{equation}
f(N)\psi_0=0,\quad f(N)\psi_1=\sqrt{2}b_0^2\psi_1,\quad
f(N)\psi_n=\sqrt{2}b_0^2\frac{v_n-v_{n-2}}{v_1}\psi_n,
\label{triv35}
\end{equation}
where $n\geq2$. The position operator  $X_{\mu}$ is defined by the
recurrence relations (\ref{triv15}) and (\ref{triv16}). Using
$X_{\mu}$ and $a_{\mu}^{-}$, we determine by the well-known
formulas (see \cite{bor}) the operators $a_{\mu}^{+}$, $P_{\mu}$
(the momentum operator) and $H_{\mu}$ (hamiltonian):
\begin{align}
\label{triv36} P_{\mu}=\imath(\sqrt{2}a_{\mu}^{-}-X_{\mu}),\qquad
a_{\mu}^{+} &=\sqrt{2}X_{\mu}-a_{\mu}^{-},\\ \label{triv37}
H_{\mu}=X_{\mu}^2+P_{\mu}^2=\sqrt{2}(a_{\mu}^{-}X_{\mu}+X_{\mu}a_{\mu}^{+})
&=a_{\mu}^{-}a_{\mu}^{+}+a_{\mu}^{+}a_{\mu}^{-}.
\end{align}
We have the following commutation relation:
\begin{equation}
[a_{\mu}^{-},a_{\mu}^{+}]=2(B(N+I)-B(N)). \label{triv38}
\end{equation}
The operator-function $f(N)$ acts on basis vectors by formulas:
\begin{equation}
B(N)\psi_0=2b_0^2,\quad
B(N)\psi_n=b_{n-1}^2\psi_n=b_0^2\frac{v_{n-1}(v_n-v_{n-2})}{v_1}\psi_n,\quad{n\geq1}.
\label{trivs39}
\end{equation}
Moreover, the "energy levels" are
\begin{equation}
\lambda_0=2b_0^2,\quad
\lambda_n=2(b_{n-1}^2+b_n^2)=\frac{2b_0^2}{v_1}(v_nv_{n+1}-v_{n-1}v_{n-2}),
\label{triv39}
\end{equation}
where $n\geq1$.

In what follows our prime interest is with the following question.
 Is it possible to get a differential equation of the second order for
 Hermite-Chihara polynomials from the equation
$H_{\mu}{\psi_n}=\lambda_n{\psi_n}$.

\section{Classical Hermite-Chihara polynomials}

Now we consider a particular case of the Hermite-Chihara
polynomials, namely, the well-known  (see \cite{sze}) generalized
Hermite polynomials which have been studied extensively in
(\cite{Chih}) (see also (\cite{det})).

We denote by ${\tt H}_{\gamma}$ the Hilbert space
\begin{equation}
{\tt
H}_{\gamma}={L^2}(R;|x|^{\gamma}(\Gamma(\frac{1}{2}(\gamma+1)))^{-1}\exp(-x^2){dx}),\qquad\gamma\geq{-1}.
\label{triv40}
\end{equation}
Using methods of (\cite{bor}), we construct the canonical
orthonormal polynomials system $\left\{
{\psi_{n}(x)}\right\}_{n=0}^{\infty }$ completed in the space
${\tt H}_{\gamma}$. The polynomials $\psi_{n}(x)$ satisfy the
recurrent relations (\ref{triv15}) and (\ref{triv16}). The
coefficients $\left\{ {b_{n}}\right\}_{n=0}^{\infty }$ are defined
by formulas (\ref{triv30}), where
 $b_0=\sqrt{\frac{\gamma+1}{2}}$ and the sequence $\vec{v}=\left\{ {v_{n}}\right\}_{n=0}^{\infty
}$ is given by the following equalities:
\begin{equation}
v_n= \Biggl\{
\begin{array}{cc}
\frac{\gamma+n+1}{\gamma+1}&{n=2m},\\
\frac{n+1}{\gamma+1}&{n=2m+1}.\\
\end{array}
\Biggr. \label{triv41}
\end{equation}
It is clear that $v_0=1,\quad{v_1=\frac{2}{\gamma+1}=b_0^{-2}}.$
The coefficients $\left\{ {b_{n}}\right\}_{n=0}^{\infty }$ are
defined by the formulas:
\begin{equation}
b_{n-1}=\frac{1}{2} \Biggl\{
\begin{array}{cc}
\sqrt{n}&{n=2m},\\ \sqrt{n+\gamma}&{n=2m+1}.\\
\end{array}
\Biggr. \label{triv42}
\end{equation}
The formulas (\ref{triv22}),(\ref{triv19}),(\ref{triv18}) give a
explicit form of the polynomials $\psi_{n}(x)$. Recall that the
polynomials
$$ K_n^{\gamma}(x)=s_n{} \psi_{n}(x),\quad {n\geq0},
$$
 as $s_0=1$ and ${s_n=(b_{n-1})!}$, are named the generalized Hermite
polynomials in (\cite{Chih})( see also (\cite{det})). In what
follows we shall call the polynomials $\psi_{n}(x)$ as the
"classical Hermite-Chihara polynomials". It is easy to prove that
the family of these polynomials is a particular case of the more
general class of Hermite-Chihara polynomials and that the
generalized derivation operator $D_{\vec{v}}$ corresponding to the
given sequence $\vec{v}$ is determined by formulas (\ref{triv4})
and (\ref{triv5}), where
\begin{equation}
\varepsilon_1=1, \qquad
\varepsilon_m={\frac{(-2)^{m-1}}{m!}}{\frac{\gamma}{\gamma+1}} ,
\qquad m\geq2 . \label{triv43}
\end{equation}
In addition, the sequence $\vec{\gamma}=\left\{
{\gamma_{n}}\right\}_{n=0}^{\infty }$ appearing in (\ref{triv17})
is defined by equalities:
\begin{equation}
\gamma_{n}=\frac{\sqrt{2}}{\gamma+1} \Biggl\{
\begin{array}{cc}
\sqrt{n}&{n=2m},\\ \sqrt{n+\gamma}&{n=2m+1}.\\
\end{array}
\Biggr. \label{triv44}
\end{equation}
Comparing (\ref{triv34}), (\ref{triv35}) and (\ref{triv17}), we
obtain
\begin{equation}
a_{\mu}^{-}=\frac{\gamma+1}{\sqrt{2}}D_{\vec{v}}. \label{triv45}
\end{equation}
The following formulas are known
(\cite{Chih})( see, also(\cite{det})):
\begin{align}
\label{triv46} \frac{d}{dx}\psi_0=0,\qquad \frac{d}{dx}\psi_n
&=\frac{n}{b_{n-1}}\psi_{n-1}+\frac{(n-1)\theta_n}{2b_{n-1}b_{n-2}}X^{-1}\psi_{n-2}=\\
\label{triv47}
 &=2b_{n-2}\psi_{n-1}-\frac{\theta_n}{x}\psi_n,\qquad{n\geq1},
\end{align}
where
\begin{equation}
\theta_n={\theta_n}(\gamma)=\gamma\frac{1-(-1)^n}{2}.
\label{triv48}
\end{equation}
Taking into account how the annihilation operator $a_{\mu}^{-}$
and the number operator $N$ act on the basis vectors $\left\{
{\psi_{n}(x)}\right\}_{n=0}^{\infty }$, it is easy to get from
 the relations (\ref{triv46})-(\ref{triv48}) the
following formula:
\begin{equation}
X_{\mu}\frac{d}{dx}-N=(a_{\mu}^{-})^{2}. \label{triv49}
\end{equation}
Note also that the action of the position operator $X_{\mu}$ on
the basis vectors $\left\{ {\psi_{n}(x)}\right\}_{n=0}^{\infty }$
is defined by (\ref{triv15}) and (\ref{triv42}). Now we consider
the operator
\begin{equation}
\Theta_N=2B(N)-N, \label{triv50}
\end{equation}
where the operator-function $B(N)$
is defined by (\ref{triv39}).
Using (\ref{triv48}) and (\ref{triv50}), we see that
\begin{equation}
\Theta_N{\psi_n}=\theta_n{\psi_n},\qquad{n\geq1}. \label{triv51}
\end{equation}
Taking into account the equation
$H_{\mu}{\psi_n}=\lambda_n{\psi_n}$, where a hamiltonian $H_{\mu}$
defined by (\ref{triv37}), and the equality
$a_{\mu}^{+}a_{\mu}^{-}=2B(N)$, we have the following relation:
\begin{equation}
a_{\mu}^{-}a_{\mu}^{+}=2B(N+I). \label{triv52}
\end{equation}
Moreover, from (\ref{triv47}) it follows that
\begin{equation}
a_{\mu}^{-}=\frac{1}{\sqrt{2}}\frac{d}{dx}+
\frac{1}{\sqrt{2}}X_{\mu}^{-1}\Theta_N.
\label{triv53}
\end{equation}
Now from (\ref{triv52}), (\ref{triv53}) and  (\ref{triv37}) we
have
\begin{equation}
(\frac{d}{dx}+X_{\mu}^{-1}\Theta_N)(X_{\mu}-
\frac{1}{2}\frac{d}{dx}-\frac{1}{2}X_{\mu}^{-1}\Theta_N)=2B(N+I).
\label{triv54}
\end{equation}
Multiplying both sides of (\ref{triv54}) by $-2X_{\mu}$ from the
left, we obtain:
\begin{align}
 -2X_{\mu}(I+X_{\mu}\frac{d}{dx})&+X_{\mu}\frac{d^2}{dx^2}
+X_{\mu}(-X_{\mu}^{-2}\Theta_N+X_{\mu}^{-1}\Theta_N)\\
-2\Theta_N{X_{\mu}}+\Theta_N{\frac{d}{dx}}&+\Theta_N{X_{\mu}^{-1}}
\Theta_N=-2X_{\mu}2B(N+I).
\label{triv55}
\end{align}
It is not hard to prove that :
\begin{align}
\label{triv56} \Theta_N{X_{\mu}^{-1}}\Theta_N\psi_n=0,\qquad
\Theta_N{X_{\mu}} &=X_{\mu}\Theta_{N+I},\qquad\\ \label{triv57}
2X_{\mu}(2B(N+I)-\Theta_{N+I})-2X_{\mu} &=2X_{\mu}N,\qquad\\
\label{triv58}
(\Theta_N\frac{d}{dx}+\frac{d}{dx}\Theta_N)\psi_n&
=\gamma\psi_n^{\prime}.
\end{align}
Applying  both sides of (\ref{triv55}) to $\psi_n$ and using
(\ref{triv56})-(\ref{triv58}), we get the following differential
equation
\begin{equation}
x\psi_n^{\prime\prime}+(\gamma-2x^2)\psi_n^{\prime}+
(2nx-\frac{\theta_n}{x})\psi_n=0,\qquad{n\geq0},
\label{triv59}
\end{equation}
which is coincident with the well-known differential equation for
the classical Hermite-Chihara polynomials \cite{Chih} (see also
\cite{sze}).
\begin{rem}
The generators $a_{\mu}^{+},a_{\mu}^{-}$ of the generalized
Heisenberg algebra $A_{\mu}$ corresponding to the classical
Hermite-Chihara polynomials system subject to the following
commutative relation (see \cite{bor}):
\begin{equation}
[a_{\mu}^{-},a_{\mu}^{+}]=(\gamma+1)I-2\Theta_N. \label{triv60}
\end{equation}
The "energy levels" of the associated oscillator are equal to:
\begin{equation}
\lambda_0=\gamma+1,\qquad \lambda_n=2n+\gamma+1,\qquad {n\geq1}.
\label{triv61}
\end{equation}
Finally, it follows from (\ref{triv37}) and (\ref{triv53}) that
the momentum operator take the following form:
\begin{equation}
P_{\mu}=\imath(\frac{d}{dx}+X_{\mu}^{-1}\Theta_N-X_{\mu}).
\label{triv62}
\end{equation}
\end{rem}

\section{Construction of the "governing sequence" $\vec{v}$ of
a generalized derivation operator $D_{\vec{v}}$ for the classical
Hermite-Chihara polynomials}

Let $\left\{ {\psi_{n}(x)}\right\}_{n=0}^{\infty }$ be a
orthonormal Hermite-Chihara polynomials system. According to
theorem {\ref{a}}, there is a sequence $\vec{v}$ such that
(\ref{triv1}) and (\ref{triv21}) are hold. The formula
(\ref{triv23})allows us to define the sequence $\vec{\gamma}$  by
 $\vec{v}$. Furthermore, the  generalized
derivation operator $D_{\vec{v}}$ corresponding to $\vec{v}$ is a
reducing operator for the system $\left\{
{\psi_{n}(x)}\right\}_{n=0}^{\infty
 }$, i.e. the equalities (\ref{triv17}) are valid. In this section
 we obtain the exact condition on $\vec{v}$ which select some
 family of  Hermite-Chihara polynomials. This family is a
natural extension of the set of classical Hermite-Chihara
polynomials. The associated set of $\vec{v}$ is a three-parameter
family depending on the parameters $b_0,v_1$ and $v_2$. But it turn
out that the parameter $v_1$ is unessential, so that  the
above-mentioned family is really a two-parameter one. According to
(\ref{triv4}) and (\ref{triv5})), the generalized derivation
operator $D_{\vec{v}}$ corresponding to a sequence $\vec{v}$
complying with (\ref{triv1}) and (\ref{triv21}) take the following
form:
\begin{equation}
D_{\vec{v}}=X^{-1}(\underline{B_1}+\overline{B_1})=\frac
{d}{dx}+X^{-1}\overline{B_1},\qquad \underline{B_1}=X\frac{d}{dx},
\label{triv63}
\end{equation}
where
\begin{equation}
\overline{B_1}=\sum_{k=2}^{\infty}{\varepsilon_k{x^k}\frac{d^k}{dx^k}}.
\label{triv64}
\end{equation}
The coefficients $\varepsilon_k$ in (\ref{triv64}) are defined
from given sequence $\vec{v}$ by the recurrent relations
(\ref{triv5}).
\begin{rem}
For the classical Hermite-Chihara polynomials it follows from
(\ref{triv64}) and (\ref{triv43}) that
\begin{equation}
X\frac{d}{dx}+{\overline{B_1}}^{cl}=\delta(N),
 \label{triv65}
\end{equation}
where $\delta(N)$ is the projection on the subspace of the
polynomials of odd degree: $$ \delta(N)x^n={\theta_n}(1)x^n,\quad
{n\geq0},
 $$
(see (\ref{triv48})) and hence
$$
\delta(N)\psi_n^{cl}={\theta_n}(1)\psi_n^{cl}.
$$
\end{rem}
We denote by $\left\{ {\psi_{n}(x)}\right\}_{n=0}^{\infty }$ the
 orthonormal Hermite-Chihara polynomials system which is
 constructed according to remark 3.4 for given sequence $\vec{v}$.
Obviously, $$ \overline{B_1}\psi_0=\overline{B_1}\psi_1=0. $$

 Now we shall restrict our consideration to a particular class of
 the Hermite-Chihara polynomials for which there are two real sequence
$\vec{\delta}=\left\{
{\overline{\delta_{n}}}\right\}_{n=2}^{\infty }$ and
$\vec{\beta}=\left\{ {\overline{\beta_{n}}}\right\}_{n=0}^{\infty}$
such that:
\begin{equation}
\overline{B_1}\psi_2=\overline{\delta_{2}}X\psi_1,\qquad
\overline{B_1}\psi_n=\overline{\delta_{n}}X\psi_{n-1}+
\overline{\beta_{n}}\psi_{n-2},\qquad{n\geq3}.
\label{triv66}
\end{equation}
Replacing $\overline{\delta_{n}}$ by $\underline{\delta_{n}}$ and
$\overline{\beta_{n}}$ by  $\underline{\beta_{n}}$, we see from
(\ref{triv63}) and (\ref{triv17}) that the condition
(\ref{triv66}) is valid for $\underline{B_1}$ too. Then it is
clear that the assumption (\ref{triv66}) takes the place of the
rule of derivation (\ref{triv47}). Substituting (\ref{triv64})
into (\ref{triv66}), and taking into account (\ref{triv19}), we
obtain the following relation:
\begin{align}
\sum_{k=2}^{n}{{\sum_{m=0}^{\epsilon(\frac{n-k}{2})}}
\frac{(-1)^{m}}{\sqrt{[n]!}}b_{0}^{2m-n}
{\alpha_{2m-1,n-1}}x^{n-2m}\frac{(n-m)!}{(n-2m-k)!}}&=\nonumber\\
=\overline{\beta_n}{\sum_{m=0}^{\epsilon(\frac{n-2}{2})}}
{\frac{(-1)^{m}}{\sqrt{[n-2]!}}
b_{0}^{2m-n+2} {\alpha_{2m-1,n-3}}x^{n-2m-2}}&+\nonumber\\
+\overline{\delta_n}{\sum_{m=0}^{\epsilon(\frac{n-1}{2})}}
{\frac{(-1)^{m}}{\sqrt{[n-1]!}}
b_{0}^{2m-n+1} {\alpha_{2m-1,n-2}}x^{n-2m}}&,\quad{n\geq2}.
\label{triv67}
\end{align}
Equating the coefficients at $x^n$ in the both sides of
(\ref{triv67}), we get
\begin{equation}
\overline{\delta_n}b_{n-1}={A_n}(2),\qquad{n\geq2}, \label{triv68}
\end{equation}
where we used the notation
\begin{equation}
{A_s}(m)=s!\sum_{k=m}^{s}{\frac{\varepsilon_k}{(s-k)!}},\qquad{s\geq{m}},
\label{triv69}
\end{equation}
and the coefficients  $\varepsilon_k$ are defined by formulas
(\ref{triv5}). For the classical Hermite-Chihara polynomials from
(\ref{triv43}), (\ref{triv44}) and binomial formula it follows that
(as $n\geq2$)
\begin{equation}
{\overline{\delta_{n}}}^{cl}=\frac{\sqrt{2}\gamma}{\gamma+1}
\Biggl\{
\begin{array}{cc}
\sqrt{n}&{n=2m},\\ \frac{n-1}{\sqrt{n+\gamma}}&{n=2m+1}.\\
\end{array}
\Biggr. \label{triv70}
\end{equation}
In order that to find the quantities $\beta_n$ in (\ref{triv66})
we equal the coefficients at $x^t$ in the both sides of
(\ref{triv67}) (as $0\leq{t}<n$). Obviously, a coefficients at
$x^t$ only distinct from zero when it is valid the following
condition:
\begin{equation}
n-t=2p,\qquad{1\leq{p}\leq\epsilon(\frac{n}{2})}. \label{triv71}
\end{equation}
We consider separately three cases $t=0,1,2.$

1. Let $t=0,\quad{n=2p}.$ We have $\beta_{2p}=0,\quad{p\geq1}.$

2. Let $t=1,\quad{n=2p+1}.$ We have
\begin{equation}
\overline{\delta_{2p+1}}=\frac{\sqrt{[2p]}{\alpha_{2p-3,2p-2}}}{b_0{\alpha_{2p-1,2p-1}}}\overline{\beta_{2p+1}},\qquad{p\geq1}.
\label{triv72}
\end{equation}

3. Let $t=0,\quad{n-2=2p}$. Regarding $p\geq1$, we have
\begin{align}
\alpha_{2p-1,2p+1}A_2(2)=-\overline{\beta}_{2p+2}\sqrt{[2p+2][2p+1]}{\alpha_{2p-3,2p-1}}+\nonumber\\
+\alpha_{2p-1,2p}A_{2p+2}(2). \label{triv73}
\end{align}
Taking into account that $\overline{\beta}_{2p}=0,$ we have from
here
\begin{equation}
\alpha_{2p-1,2p}A_{2p+2}(2)-\alpha_{2p-1,2p+1}A_2(2)=0,\qquad{p\geq1}.
\label{triv74}
\end{equation}
Using (\ref{triv22}), we simplify this relation:
\begin{equation}
v_1A_{2p+2}(2)-v_{2p+1}A_2(2)=0,\qquad{p\geq1}. \label{triv75}
\end{equation}
From the equalities (\ref{triv5}) and the designation
(\ref{triv69}) it follows that:
\begin{equation}
A_{k}(2)=v_{k-1}-k,\qquad{k\geq2}. \label{triv76}
\end{equation}
Substituting (\ref{triv76}) into (\ref{triv75}), we get
\begin{equation}
v_{2p+1}=(p+1)v_1,\qquad{p\geq1}. \label{triv77}
\end{equation}

Now we consider the general case $t\geq3,\quad{n\geq5}$ and
$1\leq{p}\leq{\epsilon(\frac{n-3}{2})}.$ From (\ref{triv67})-
(\ref{triv69}) we have:
\begin{equation}
\alpha_{2p-1,n-1}A_{n-2p}(2)=
-\overline{\beta}_{n}\sqrt{[n][n-1]}{\alpha_{2p-3,n-3}}+
\alpha_{2p-1,n-2}A_{n}(2).
\label{triv78}
\end{equation}
As $n\geq2$ and ${2p}\leq{n}$, from the condition (3.5) of the
paper (\cite{bodak}) it follows that:
\begin{equation}
\alpha_{2p-1,n-1}=[n-1]{\alpha_{2p-3,n-3}}+\alpha_{2p-1,n-2}.
\label{triv79}
\end{equation}
Substituting (\ref{triv79}) into (\ref{triv78}), we get the
following formula:
\begin{align}
\alpha_{2p-3,n-3}([n-1]A_{n-2p}(2)+\overline{\beta}_{n}
\sqrt{[n][n-1]})=\nonumber\\
=\alpha_{2p-1,n-2}(A_{n}(2)-A_{n-2p}(2)). \label{triv80}
\end{align}
Combining (\ref{triv22}),(\ref{triv29}) and (\ref{triv30}),
we get (as $p\geq1$ and ${n\geq{2p+1}}$)
\begin{equation}
\frac{\alpha_{2p-1,n-2}}{\alpha_{2p-3,n-3}}=
\frac{v_{2p-1}-v_{2p-3}}{v_1v_{2p-1}}v_{n-2}v_{n-2p-1},
\label{triv81}
\end{equation}
as well as
\begin{equation}
[n-1]=\frac{v_{n-2}(v_{n-1}-v_{n-3})}{v_1}. \label{triv82}
\end{equation}
Substituting (\ref{triv81}) into (\ref{triv80}), we get (as
$n\geq5$ and ${1\leq{p}\leq{\epsilon(\frac{n-3}{2})}}$) the
following formula:
\begin{align}
-\frac{\overline{\beta}_{n}\sqrt{[n][n-1]}}{v_{n-2}}=
\frac{v_{n-1}-v_{n-3}}{v_1}A_{n-2p}(2)-\nonumber\\
-\frac
{v_{n-2p}(v_{2p-1}-v_{2p-3})}{v_1v_{2p-1}}(A_{n}(2)-A_{n-2p}(2)) .
\label{triv83}
\end{align}
Taking into account (\ref{triv21}), the coefficient at
$A_{n-2p}(2)$ in the right side of (\ref{triv83}) is equal to
\begin{align}
\frac{(v_{n-1}-v_{n-3})v_{2p-1}+v_{n-2p-1}(v_{2p-1}-v_{2p-3})}
{v_1v_{2p-1}}\nonumber\\
=\frac{(v_{2p-1}-v_{2p-3})v_{n-1}}{v_1v_{2p-1}} . \label{triv84}
\end{align}
As $n\geq5$ and ${1\leq{p}\leq{\epsilon(\frac{n-3}{2})}}$, from
(\ref{triv83}) and (\ref{triv84}) we get
\begin{equation}
-\frac{\overline{\beta}_{n}\sqrt{[n][n-1]}}{v_{n-2}}=
\frac{v_{2p-1}-v_{2p-3}}{v_1v_{2p-1}}
(v_{n-1}A_{n-2p}(2)- A_{n}(2)v_{n-2p-1}) . \label{triv85}
\end{equation}
As $n\geq5$ and ${1\leq{p}\leq{\epsilon(\frac{n-3}{2})}}$, taking
into account (\ref{triv76}) and (\ref{triv77}), we can rewrite the
right hand side of (\ref{triv85}) in the following form:
\begin{equation}
-\frac{\overline{\beta}_{n}\sqrt{[n][n-1]}}{v_{n-2}}=
\frac{nv_{n-2p-1}-(n-2p)v_{n-1}}{pv_1}.
\label{triv86}
\end{equation}
Evidently, the condition (\ref{triv86}) is valid as $n=2p$, if it
is remembered that $\beta_{2p}=0$ as well as (\ref{triv77}). It
remains to consider the case $n=2m+1$ ($m\geq2$ and
${p\leq{m-1}}$):
\begin{equation}
-\frac{\overline{\beta}_{2m+1}\sqrt{[2m+1][2m]}}{v_{2m-1}}=
\frac{(2m+1)v_{2m-2p}-(2m-2p+1)v_{2m}}{pv_1}.
\label{triv87}
\end{equation}
From (\ref{triv72}), (\ref{triv76}) and (\ref{triv81}), as
$n=2m+1$ ,$m\geq2$ and ${p\leq{m-1}}$, we have
\begin{align}
\overline{\beta}_{2m+1}{\sqrt{[2m]}}
={b_0}\overline{\delta}_{2m+1}
\frac{\alpha_{2m-1,2m-1}}{\alpha_{2m-3,2m-2}}=\nonumber\\
={b_0}\overline{\delta}_{2m+1}=\frac{A_{2m+1}}{\sqrt{[2m+1]}}=
\frac{v_{2m}-(2m+1)}{\sqrt{[2m+1]}}.
\label{triv88}
\end{align}
Substituting (\ref{triv88}) into (\ref{triv87}) and using
(\ref{triv77}), we obtain the following relation:
\begin{equation}
-(v_{2m}-(2m+1))=\frac{m}{p}((2m+1)v_2(m-p)-(2m-2p+1)v_{2m}),
\label{triv89}
\end{equation}
i.e. for any $m\geq2$ and ${p\leq{m-1}}$ it should be true the
following equality:
\begin{equation}
(m-p)v_{2m}+p=mv_2(m-p). \label{triv90}
\end{equation}
It is readily seen that (\ref{triv90}) is valid if and only if:
\begin{equation}
v_{2m}=mv_2-(m-1),\qquad {m\geq1}. \label{triv91}
\end{equation}
So it is proved the following theorem.
\begin{thm}\label{e}
For the orthonormal  Hermite-Chihara  polynomials system
constructed by the sequence $\vec{v}$, submitting to (\ref{triv1})
and (\ref{triv21}), be satisfied (\ref{triv66}) it is necessary
and sufficient that this system is obeying also the relations
(\ref{triv77}) and (\ref{triv91}) at some $v_1$ and $v_2$ such
that $1\leq{v_1}\leq{v_2}$.
\end{thm}
\begin{rem}
Thus we constructed the three-parameter Hermite-Chihara
polynomials family (these parameters are $b_0,v_1$ and $v_2$)
 for which is hold the condition (\ref{triv66}). For the classical
Hermite-Chihara polynomials are valid the following relations:
\begin{equation}
v_1={b_0}^{-2}=\frac{2}{\gamma+1},\qquad v_2=1+v_1. \label{triv92}
\end{equation}
\end{rem}

\section {Deduction of the differential equation for
the  family of Hermite-Chihara polynomials}

In this section we prove that any polynomial belonging to the
three-parameter family considered above satisfies a
the second order differential equation.

From (\ref{triv17}), (\ref{triv63}) and (\ref{triv66}) it follows
that for any polynomials system $\left\{
{\psi_{n}(x)}\right\}_{n=0}^{\infty }$ belonging to the considered
family the differentiation operator $\frac{d}{dx}$ acts on this
basis by the folowing formulas:
\begin{align}
\frac{d}{dx}\psi_0=0,\qquad
\frac{d}{dx}\psi_1=\frac{1}{b_0},\qquad\nonumber\\
\frac{d}{dx}\psi_n=(\gamma_n-
\overline{\delta}_{n})\psi_{n-1}-\overline{\beta}_{n}
X_{\mu}^{-1}\psi_{n-2},\quad{n\geq2}.
\label{triv93}
\end{align}
Multiplying (\ref{triv93}) from the left by $X_{\mu}$  and using
(\ref{triv15}), we get
\begin{equation}
(X_{\mu}\frac{d}{dx}-b_{n-1}(\gamma_n-
\overline{\delta}_{n}))\psi_n=(b_{n-2}(\gamma_n-
\overline{\delta}_{n})-\overline{\beta}_{n})\psi_{n-2},\quad
{n\geq2}. \label{triv94}
\end{equation}
Note that from (\ref{triv23}) and (\ref{triv30}) it follows the
relation:
\begin{equation}
v_{n-1}=b_{n-1}\gamma_n,\qquad{n\geq1}. \label{triv95}
\end{equation}
Then from (\ref{triv95}), (\ref{triv68}) and (\ref{triv76}) we
have as $n\geq2$
\begin{equation}
b_{n-1}(\gamma_n- \overline{\delta}_{n})=v_{n-1}-A_n(2)=n,
\label{triv96}
\end{equation}
and hence
\begin{equation}
b_{n-2}(\gamma_n-
\overline{\delta}_{n})=n\frac{b_{n-2}}{b_{n-1}},\qquad{n\geq2}.
\label{triv97}
\end{equation}
From (\ref{triv93}),(\ref{triv94}), (\ref{triv96}) and
(\ref{triv97}) we have
\begin{align}
\frac{d}{dx}\psi_0=0,\qquad
\frac{d}{dx}\psi_1=\frac{1}{b_0},\qquad\nonumber\\
(X_{\mu}\frac{d}{dx}-N)\psi_n=(-\overline{\beta}_{n}+
n\frac{b_{n-2}}{b_{n-1}})\psi_{n-2},\quad{n\geq2}.
\label{triv98}
\end{align}
Let us remark that from (\ref{triv68}) and (\ref{triv72})
\begin{equation}
\overline{\beta}_{2p+1}=b_0\frac{A_{2p+1}(2)\alpha_{2p-1,2p-1}}
{b_{2p}\sqrt{[2p]}\alpha_{2p-3,2p-2}}.
\label{triv99}
\end{equation}
By virtue of (\ref{triv81}) as $n=2p+1$ and taking into account
(\ref{triv77}), we have
\begin{equation}
\frac{\alpha_{2p-1,2p-1}}{\alpha_{2p-3,2p-2}}=\frac{v_{2p-1}-
v_{2p-3}}{v_1v_{2p-1}}v_{2p-1}=1.
\label{trivs1}
\end{equation}
Further, from (\ref{triv99}) , (\ref{trivs1}),
(\ref{triv76}),(\ref{triv91}) and (\ref{triv18}) it follows that
\begin{equation}
\overline{\beta}_{2p+1}={b_0}^2\frac{v_{2p}-(2p+1)}{b_{2p}b_{2p-1}}
={b_0}^2\frac{(v_{2}-3)p}{b_{2p}b_{2p-1}}.
\label{trivs2}
\end{equation}
Note also that from (\ref{triv30}),(\ref{triv91}) and
(\ref{triv77}) it follows that
\begin{align}
b_{2p-1}^2={b_0}^2\frac{v_{2p-1}(v_{2p}-v_{2p-2})}{v_1}={b_0}^2(v_2-1)p,
\qquad\nonumber\\
b_{2p}^2={b_0}^2\frac{v_{2p}(v_{2p+1}-v_{2p-1})}{v_1}={b_0}^2v_{2p}=
{b_0}^2(pv_2-(p-1)).
\label{trivs3}
\end{align}
Taking into account that $\beta_{2p}=0$, the right side of
  (\ref{triv98}) takes the form:
\begin{equation}
-\overline{\beta}_{n}+n\frac{b_{n-2}}{b_{n-1}}= \Biggl\{
\begin{array}{cc}
\frac{nb_{n-2}}{b_{n-1}}&{n=2m},\\
\frac{(n-1)b_{n-1}}{b_{n-2}}&{n=2m+1}.\\
\end{array}
\Biggr. \label{trivs4}
\end{equation}
Then as $n=2p$
\begin{equation}
\frac{{2p}b_{2p-2}}{b_{2p-1}}=\frac{{2p}b_{2p-2}b_{2p-1}}{b_{2p-1}^2}
=\frac{2b_{n-1}b_{n-2}}{b_0^2(v_2-1)},
\label{trivs5}
\end{equation}
and as $n=2p+1$
\begin{equation}
\frac{{2p}b_{2p}}{b_{2p-1}}=\frac{{2p}b_{2p-2}b_{2p-1}}{b_{2p-1}^2}=
\frac{2b_{n-1}b_{n-2}}{b_0^2(v_2-1)},
\label{trivs6}
\end{equation}
so from (\ref{trivs5}),(\ref{trivs6}) and (\ref{triv98}) it
follows that
\begin{equation}
X_{\mu}\frac{d}{dx}-N=c_1^{-1}(a_{\mu}^{-})^2, \label{trivs7}
\end{equation}
where
\begin{equation}
c_1=b_0^2(v_2-1). \label{trivs8}
\end{equation}
We stress that, according to (\ref{triv92}), the formula
(\ref{trivs7}) is a extension of the one (\ref{triv49}).

Next we will arguing by analogy with the derivation of the
differential equation (\ref{triv59}). According to (\ref{triv37}),
(\ref{triv52}) and using (\ref{trivs7}) as well as the following
notation:
\begin{equation}
\Delta_N=2{c_1^{-1}}B(N)-N, \label{trivs9}
\end{equation}
we get:
\begin{equation}
a_{\mu}^{-}a_{\mu}^{+}=c_1(\frac{d}{dx}+X_{\mu}^{-1}\Delta_N)
(X_{\mu}-\frac{c_1}{2}(\frac{d}{dx}+X_{\mu}^{-1}\Delta_N))=2B(N+I).
\label{trivs10}
\end{equation}
Taking into account the notation (\ref{trivs8}), we rewrite
(\ref{trivs3}):
\begin{equation}
2{b_{n-1}}^2= \Biggl\{
\begin{array}{cc}
nc_1,&{n=2m},\\ nc_1-c_1+2{b_0}^2,&{n=2m+1}.\\
\end{array}
\Biggr. \label{trivs11}
\end{equation}
From (\ref{trivs9}),(\ref{trivs39}), (\ref{triv96}) and
(\ref{trivs11}) we have:
\begin{equation}
\Delta_N^{c_1}=\alpha_n\psi_n, \label{trivs12}
\end{equation}
where
\begin{equation}
\alpha_n=\frac{3-v_2}{v_2-1}\frac{1-(-1)^n}{2}. \label{trivs13}
\end{equation}
Then can one rewrite the equation (\ref{trivs10}) in the form:
\begin{equation}
(\frac{d}{dx}+X_{\mu}^{-1}\Delta_N^{c_1})(X_{\mu}-\frac{c_1}{2}
(\frac{d}{dx}+X_{\mu}^{-1}\Delta_N^{c_1}))=2{c_1}^{-1}B(N+I).
\label{trivs14}
\end{equation}
Multiplying (\ref{trivs14}) by $-\frac{2}{c_1}X_{\mu}$ from the
left, we have:
\begin{align}
-\frac{2}{c_1}X_{\mu}(I+X_{\mu}\frac{d}{dx})+X_{\mu}(-X_{\mu}^{-2}
\Delta_N^{c_1}+X_{\mu}^{-1
}\frac{d}{dx}\Delta_N^{c_1})- \qquad\nonumber\\
-\frac{2}{c_1}\Delta_N^{c_1}X_{\mu}+\Delta_N^{c_1}\frac{d}{dx}+
\Delta_N^{c_1}X_{\mu}^{-1}\Delta_N^{c_1}=-\frac{2}{c_1}X_{\mu}
\frac{2}{c_1}B(N+I).
\label{trivs15}
\end{align}
Applying (\ref{trivs15}) to $\psi_n$ and using the following
equalities:
\begin{align}
\label{trivs16} \Delta_N^{c_1}X_{\mu}^{-1 }\Delta_N^{c_1}=0,
\qquad\\ \label{trivs17}
\Delta_N^{c_1}X_{\mu}=X_{\mu}\Delta_{N+1}^{c_1},\qquad\\
\label{trivs18}
\frac{2}{c_1}X_{\mu}(\frac{2}{c_1}B(N+I)-\Delta_{N+1}^{c_1})-
\frac{2}{c_1}X_{\mu}=\frac{2}{c_1}X_{\mu}{N},\qquad\\
\label{trivs19}
(\Delta_N^{c_1}X_{\mu}^\frac{d}{dx}+\frac{d}{dx}\Delta_N^{c_1})
\psi_n=\frac{3-v_2}{v_2-1}\psi_n^{\prime},\qquad
\end{align}
we get the following differential equation (as all $n\geq0$):
\begin{equation}
x\psi_n^{\prime\prime}+(\frac{3-v_2}{v_2-1}-
\frac{2}{b_0^2(v_2-1)}x^2)\psi_n^{\prime}+(\frac{2}{b_0^2(v_2-1)}{n}x-
\frac{\alpha_n}{x})\psi_n=0.
\label{trivs20}
\end{equation}
Because  $\alpha_n^{cl}=\theta_n$ è $c_1^{cl}=1$ it is clear that
the equation (\ref{trivs15}) is a extension of the differential
equation (\ref{triv59}). Also, we obtain the following
two-parameter differential equation (as any  $n\geq0$):
\begin{equation}
x\psi_n^{\prime\prime}+(\gamma-2\alpha{x^2})\psi_n^{\prime}+
(2\alpha{n}x-\frac{\theta_n(\gamma)}{x})\psi_n=0.
\label{trivs20}
\end{equation}
where
\begin{equation}
\gamma=\frac{3-v_2}{v_2-1},\quad
\alpha=\frac{1}{b_0^2(v_2-1)},\quad{\gamma>-1},\quad{\alpha>0},
\label{trivs21}
\end{equation}
for the  Hermite-Chihara polynomials
$\left\{{\psi_{n}(x)}\right\}_{n=0}^{\infty }$ which are
orthonormal with respect to the measure
\begin{align}
d\mu(x)=C|x|^{\gamma}\exp(-\alpha {x^2})dx, \qquad\nonumber\\
C=\alpha^{\frac{\gamma+\sqrt{\gamma^2+8}}{2}}
\frac{1}{\Gamma(\frac{\gamma+1}{2})}.
\label{trivs22}
\end{align}
The polynomials $\left\{{\psi_{n}(x)}\right\}_{n=0}^{\infty }$
satisfy the recurrence relations (\ref{triv15}) with the
coefficients:
\begin{equation}
2{b_{n-1}}^2= \Biggl\{
\begin{array}{cc}
\frac{n}{\alpha},&{n=2m},\\ \frac{n+\gamma}{\alpha},&{n=2m+1}.\\
\end{array}
\Biggr. \label{trivs23}
\end{equation}
Notice that the classical  Hermite-Chihara polynomials correspond
to the case $\alpha=1.$

\section{Conclusion}

It is easily shown that the condition (\ref{triv66}) is equivalent
to the following conditions:
\begin{align}
\underline{B_1}\psi_n=\underline{a_{n}}\psi_{n}+
\underline{c_{n}}\psi_{n-2}, \qquad\nonumber\\
\overline{B_1}\psi_n=\overline{a_{n}}\psi_{n}+
\overline{c_{n}}\psi_{n-2},\qquad{n\geq2},
\label{trivs24}
\end{align}
where from (\ref{triv63}), (\ref{triv17}) and (\ref{triv15}) it
follows that:
\begin{equation}
\overline{a_{2}}=0,\quad
\overline{a_n}+\underline{a_{n}}=\gamma_n{b_{n-1}},\quad
\overline{c_n}+\underline{c_{n}}=\gamma_n{b_{n-2}},\quad{n\geq2}.
\label{trivs25}
\end{equation}
As the operator $\overline{B_1}$ and hence the operator
$\underline{B_1}$ is a reduced operator (the lowering on the basis
vectors shall be no more than two steps), so the the differential
equation (\ref{trivs20}) is a  differential equation of the second
order. Assume that (\ref{trivs24}) does not hold, i.e., we
consider the general case:
\begin{align}
\underline{B_1}\psi_n=\underline{a_{n}}\psi_{n}+\underline{c_{n}}
\psi_{n-2}+\underline{d_{n,2}}\psi_{n-4}+\cdots+
\underline{d_{n,\epsilon(\frac{n}{2})}}
\psi_{n-2{\epsilon(\frac{n}{2})}},\nonumber\\
\overline{B_1}\psi_n=\overline{a_{n}}\psi_{n}+
\overline{c_{n}}\psi_{n-2}+\overline{d_{n,2}}
\psi_{n-4}+\cdots+\overline{d_{n,\epsilon(\frac{n}{2})}}
\psi_{n-2{\epsilon(\frac{n}{2})}},
\label{trivs26}
\end{align}
where (\ref{trivs25}) still stand as before. Besides, we have
\begin{equation}
\overline{d_{n,s}}+\underline{d_{n,s}}=0,\quad
s=2,3,...\epsilon(\frac{n}{2}), \label{trivs27}
\end{equation}
with the coefficients $\overline{d_{n,s}}$ (correspondingly
$\underline{d_{n,s}}$) do not all equal to zero. One can prove
that if one of the coefficients vanish as fixed $n$,i.e.,
\begin{equation}
\underline{d_{n,s}}=0,\quad s\geq2,n\geq2, \label{trivs28}
\end{equation}
so it follows that for all $n\geq2$ the same is true:
\begin{equation}
\underline{d_{n,s}}=0. \label{trivs29}
\end{equation}
In addition, as $n=2k$ from
\begin{equation}
\underline{d_{n,2}}=\underline{d_{n,3}}=\cdots=\underline{d_{n,s}}=0,
\label{trivs30}
\end{equation}
it follows that (\ref{triv77}) is true for all  $1\leq{p}\leq{s}$.
In the case $n=2m+1$ from (\ref{trivs30}) it follows that
(\ref{triv91}) is valid for all $1\leq{m}\leq{s}$.

Furthermore, it seems that the following assumption hold. The
Hermite-Chihara polynomials either satisfy the differential
equation (\ref{trivs20}) of the second order (if the condition
(\ref{trivs24}) is hold) or they do not  satisfy  any differential
equation of a finite order (if the general condition
(\ref{trivs26}) is hold). The proof  of  the promoted conjecture
invites further investigation and this will be the object of
another paper.

In conclusion we note that the description of the generalized
Hahn-Hermite polynomials (see (\cite{hah})) in the framework of
the our scheme as thought is related either with a relaxation  of
the condition (\ref{triv17}) or (what is more radical) with a
relaxation  of the condition (\ref{triv2}). This will be discussed
elsewhere.

\bibliographystyle{amsplain}
\renewcommand{\refname}{References}

\end{document}